\documentclass{amsart}
\usepackage{amssymb}
\usepackage[all]{xy}

\CompileMatrices 

\begin{document}

\def\T{\mathbb T}
\def\N{\mathbb N}
\def\NN{\mathbb N}
\def\ZZ{\mathbb Z}
\def\Oo{\mathcal O}
\def\Cc{\mathcal C}

\def\ker{\operatorname{ker}}
\def\coker{\operatorname{coker}}
\def\rank{\operatorname{rank}}
\def\tor{\operatorname{tor}}
\def\im{\operatorname{Im}}
\def\Obj{\operatorname{Obj}} 
\def\Mor{\operatorname{Mor}}
\def\Id{\operatorname{id}}
\def\dom{\operatorname{dom}}
\def\cod{\operatorname{cod}}
\def\Aut{\operatorname{Aut}}
\def\Span{\operatorname{span}}

\newcommand{\putbox}[4]{\multiput(#1,#2)(#3,0){2}{\line(0,1){#4}}\multiput(#1,#2)(0,#4){2}{\line(1,0){#3}}}

\def\L{\Lambda}
\def\One{{\mathbf{1}}}

\theoremstyle{plain}
\newtheorem{theorem}{Theorem}[section]
\newtheorem*{theorem*}{Theorem}
\newtheorem*{prop*}{Proposition}
\newtheorem{cor}[theorem]{Corollary}
\newtheorem{lemma}[theorem]{Lemma}
\newtheorem{prop}[theorem]{Proposition}
\theoremstyle{remark}
\newtheorem{rmk}[theorem]{Remark}
\newtheorem{rmks}[theorem]{Remarks}
\newtheorem{example}[theorem]{Example}
\newtheorem{examples}[theorem]{Examples}
\theoremstyle{definition}
\newtheorem{dfn}[theorem]{Definition}
\newtheorem{dfns}[theorem]{Definitions}
\newtheorem{notation}[theorem]{Notation}

\numberwithin{equation}{section}

\title
[Higher-rank dual graphs and $K$-theory]
{A dual graph construction for higher-rank graphs, and $K$-theory for finite 2-graphs}
\author{Stephen Allen}
\author{David Pask}
\author{Aidan Sims}
\address{Department of Mathematics  \\
      University of Newcastle\\
NSW  2308\\
AUSTRALIA}
\email{stephen.allen@studentmail.newcastle.edu.au}
\email{david.pask@newcastle.edu.au}
\email{aidan@frey.newcastle.edu.au}
\date{February 6, 2004}
\thanks{This research was supported by the Australian Research Council.}

\begin{abstract} 
Given a $k$-graph $\L$ and an element $p$ of $\NN^k$, we define the dual $k$-graph, $p\L$. We show that when $\L$ is row-finite and has no sources, the $C^*$-algebras $C^*(\L)$ and $C^*(p\L)$ coincide. We use this isomorphism to apply Robertson and Steger's results to calculate the $K$-theory of $C^*(\L)$ when $\L$ is finite and strongly connected and satisfies the aperiodicity condition.
\end{abstract}

\maketitle

\section{Introduction}
In 1980, Cuntz and Krieger introduced a class of $C^*$-algebras, now called Cuntz-Krieger algebras, associated to finite $\{0,1\}$-matrices $A$ \cite{CK}. Enomoto and Watatani then showed that these algebras could be regarded as being associated in a natural way to finite directed graphs by regarding $A$ as the vertex adjacency matrix of a finite directed graph $E$ \cite{EW}. Generalising this association, Enomoto and Watatini associated $C^*$-algebras $C^*(E)$ to finite graphs $E$ with no sources\footnote{For the sake of consistency with $k$-graph notation, we regard directed graphs as $1$-graphs, so \emph{no sources} here corresponds to \emph{no sinks} in, for example, \cite{EW,BPRS}} ($E$ has no sources if each vertex of $E$ is the range of at least one edge). Although not every finite directed graph with no sources has a vertex adjacency matrix with entries in $\{0,1\}$, the vertex adjacency matrix of the dual graph $\widehat E$ formed by regarding the edges of $E$ as vertices and the paths of length 2 in $E$ as edges \emph{does} always have entries in $\{0,1\}$, and the Cuntz-Krieger algebras associated to $E$ and to $\widehat E$ are canonically isomorphic \cite{MRS}. These results have since been extended to infinite graphs (see for example \cite{KPRR, KPR, BPRS, FLR}; see also \cite{BP} when $E$ has sources).

One of the major attractions of graph algebras is their applicability to the classification program for simple purely infinite nuclear $C^*$-algebras. Conditions on a graph $E$ have been identified which guarantee that $C^*(E)$ is  purely infinite, simple, and nuclear, and satisfies the Universal Coefficient Theorem (see, for example, \cite{BPRS}), thus producing a large class of directed graphs whose $C^*$-algebras are determined up to isomorphism by their $K$-theory \cite{P}. The $K$-theory of $C^*(E)$ for an arbitrary directed graph $E$ was calculated in \cite{RaSz}, and it is shown in \cite{Sz} that given any two finitely generated abelian groups $G, H$ such that $H$ is torsion-free, there exists a directed graph $E$ such that $C^*(E)$ is simple, purely infinite, nuclear, and satisfies the Universal Coefficient Theorem, with  $K_0(C^*(E)) \cong G$ and $K_1(C^*(E)) \cong H$.

In 1999, Robertson and Steger introduced a class of higher-rank Cuntz-Krieger algebras $\mathcal{A}$, associated to collections $M_1,\dots,M_k$ of commuting $\{0,1\}$-matrices satisfying appropriate compatibility conditions \cite{RS'}. In \cite{RS}, they went on to calculate the $K$-theory of $\mathcal{A}$, demonstrating in particular that $K_1(A)$ need not be torsion-free, so that the class of higher-rank Cuntz-Krieger algebras exhausts some $K$-invariants which are not achieved by graph algebras. In order to place these higher-rank Cuntz-Krieger algebras in a graph-theoretic setting, and to generalise them as Watatani and Enomoto had generalised the original Cuntz-Krieger algebras, Kumjian and Pask introduced the notion of a higher-rank graph $\L$, and defined and investigated the associated higher-rank graph $C^*$-algebra $C^*(\L)$ \cite{KP}. Connectivity in a rank-$k$ graph $\L$ is described in terms of $k$ commuting vertex adjacency matrices $\{M^\L_1, \dots, M^\L_k\}$, called coordinate matrices. Just as in the rank-1 setting, not every $k$-graph has coordinate matrices with entries in $\{0,1\}$, but if $\L$ is a $k$-graph whose coordinate matrices are $\{0,1\}$-matrices, then \cite[Corollary~3.5(ii)]{KP} shows that $C^*(\L)$ and the $C^*$-algebra $\mathcal{A}$ associated to the coordinate matrices as in \cite{RS'} are identical.

In this note we introduce a notion of a dual graph for higher-rank graphs, and show that for a large class of higher-rank graphs $\L$, the dual higher-rank graph $p\L$ and the original higher-rank graph $\L$ have canonically isomorphic $C^*$-algebras for all $p \in \NN^k$ (c.f. \cite{B}). We also show that by choosing $p$ appropriately, we can ensure that $p\L$ has coordinate matrices with entries in $\{0,1\}$. Using these results, we identify a class of finite rank-2 graphs whose $C^*$-algebras are isomorphic to the rank-2 Cuntz-Krieger algebras studied by Robertson and Steger, and we use the results of \cite{RS} to show that these $C^*$-algebras are purely infinite, simple, unital and nuclear, and to calculate their $K$-theory.

The layout of the paper is as follows: in Section~\ref{sec:prelims}, we recall the definition of $k$-graphs and the associated notation; in Section~\ref{sec:duals}, we introduce the dual graph construction for $k$-graphs, and show that this construction preserves the associated $C^*$-algebra; and in Section~\ref{sec:K-theory}, we identify the finite 2-graphs $\L$ whose $C^*$-algebras can be studied using Robertson and Steger's results, and use these results to calculate $K_*(C^*(\L))$.

\smallskip

In the final stages of preparation of this paper, the authors became aware of Evans' Ph.D. thesis \cite{Gw}, which appears to obtain more general results regarding $K$-theory for $2$-graph $C^*$-algebras than those established here.

\section{Preliminaries}\label{sec:prelims}
We regard $\NN^k$ as an additive semigroup with identity 0. Given $m,n \in \NN^k$, we write $m \vee n$ for their coordinate-wise maximum and $m \wedge n$ for their coordinate-wise minimum, and if $m \le n$, then we write $[m,n]$ for the set $\{p \in \NN^k : m \le p \le n\}$. We denote the canonical generators of $\NN^k$ by $\{e_1,\dots,e_k\}$, and for $n \in \NN^k$, we write $n_j$ for the $j^{\rm th}$ coordinate of $n$).

\begin{dfn} \label{dfn:k-graph} 
Let $k \in \NN \setminus \{0\}$.  A \emph{$k$-graph}, is a pair $(\L,d)$ where $\L$ is a countable category and $d$ is a functor from $\L$ to $\NN^k$ which satisfies the \emph{factorisation property}:  if $\lambda \in \Mor(\L)$ and $d(\lambda) = m+n$, then there are unique morphisms $\mu \in d^{-1}(m)$ and $\nu \in d^{-1}(n)$ such that $\lambda = \mu\nu$.
\end{dfn}

We refer to elements of $\Mor(\L)$ as \emph{paths} and to elements of $\Obj(\L)$ as \emph{vertices} and we write $r$ and $s$ for the codomain and domain maps. The factorisation property allows us to identify $\Obj(\L)$ with $\{\lambda \in \Mor(\L) : d(\lambda) = 0\}$. So we write $\lambda \in \L$ in place of $\lambda \in \Mor(\L)$, and when $d(\lambda) = 0$, we regard $\lambda$ as a vertex of $\L$.

Given $\lambda \in \L$ and $E \subset\L$, we define $\lambda E :=\{\lambda\mu : \mu \in E, r(\mu) = s(\lambda)\}$ and $E\lambda := \{\mu\lambda : \mu \in E, s(\mu) = r(\lambda)\}$.  In particular if $d(v) = 0$, then $v$ is a vertex of $\L$ and $vE = \{\lambda \in E : r(\lambda) = v\}$; similarly, $Ev = \{\lambda \in \L : s(\lambda) = v\}$.  We write
$\L^n$ for the collection $\{\lambda \in \L : d(\lambda) = n\}$.

\begin{dfn}\label{dfn:structural conditions}
We say that a $k$-graph $(\L,d)$ is \emph{row-finite} if $v\L^n$ is finite for all $v \in \L^0$ and $n \in \NN^k$, and that $\L$ has \emph{no sources} if $v\L^n$ is nonempty for all $v \in \L^0$ and $n \in \NN^k$. We say that $\L$ is \emph{strongly connected\/} if $v\L w$ is nonempty for all $v,w \in \L^0$, and we say that $\L$ is $\emph{finite}$ if $\L^0$ and each $\L^{e_i}$ are finite.
\end{dfn}

The factorisation property ensures that if $l \le m \le n \in \NN^k$ and if $d(\lambda) = n$, then there exist unique paths denoted $\lambda(0, l)$, $\lambda(l, m)$ and $\lambda(m,n)$ such that $d(\lambda(0,l)) = l$, $d(\lambda(l,m)) = m - l$, and $d(\lambda(m,n)) = n-m$ and such that $\lambda = \lambda(0,l)\lambda(l,m)\lambda(m,n)$.

Given $k \in \NN\setminus\{0\}$, and $k$-graphs $(\L_1, d_1)$ and $(\L_2, d_2)$, we call a covariant functor $x : \L_1 \to \L_2$ a \emph{graph morphism} if it satisfies $d_2 \circ x = d_1$.

\begin{dfn}\label{dfn:infinite paths}
As in \cite{KP}, given $k \in \NN\setminus\{0\}$, we write $\Omega_k$ for the $k$-graph given by $\Obj(\Omega_k) = \NN^k$, $\Mor(\Omega_k) = \{(m,n) \in \NN^k \times \NN^k : m \le n\}$, $r(m,n) = m$, $s(m,n) = n$, $(m,n) \circ (n,p) = (m,p)$, and $d(m,n) = n-m$. Given a $k$-graph $\L$, an \emph{infinite path} of $\L$ is a graph morphism $x : \Omega_k \to \L$. We denote the collection of all infinite paths of $\L$ by $\L^\infty$. For $p \in \NN^k$, we write $\sigma^p : \L^\infty \to \L^\infty$ for the shift-map determined $\sigma^p(x)(n) = x(n+p)$, and we say that $x \in \L^\infty$ is \emph{aperiodic} if there do not exist $p,q \in \NN^k$ with $p \not= q$ and $\sigma^p(x) = \sigma^q(x)$.
\end{dfn}

\begin{dfn} \label{dfn:CKfam}
Let $(\L,d)$ be a row-finite $k$-graph with no sources. A Cuntz-Krieger $\L$-family is a collection $\{t_\lambda : \lambda \in \L\}$ of partial isometries satisfying
\begin{itemize}
\item[(i)] $\{t_v : v \in \L^0\}$ is a collection of mutually orthogonal projections;
\item[(ii)] $t_\lambda t_\mu = t_{\lambda\mu}$ whenever $s(\lambda) = r(\mu)$;
\item[(iii)] $t^*_\lambda t_\lambda = t_{s(\lambda)}$ for all $\lambda \in \L$; and
\item[(iv)] $t_v = \sum_{\lambda \in v\L^n} t_\lambda t^*_\lambda$ for all $v \in \L^0$ and $n \in \NN^k$.
\end{itemize}
The \emph{Cuntz-Krieger algebra} $C^*(\L)$ is the $C^*$-algebra generated by a Cuntz-Krieger $\L$-family $\{s_\lambda : \lambda \in \L\}$ which is universal in the sense that for every Cuntz-Krieger $\L$-family $\{t_\lambda : \lambda \in \L\}$ there is a unique homomorphism $\pi$ of $C^*(\L)$ satisfying $\pi(s_\lambda) = t_\lambda$ for all $\lambda \in \L$.
\end{dfn}

\section{Dual Higher Rank Graphs}\label{sec:duals}

In this section we define the higher rank analog $p\Lambda$ of the dual graph construction for directed graphs.

\begin{dfn}
Let $(\L,d)$ be a $k$-graph. Let $p\L := \{\lambda \in \L : d(\lambda) \ge p\}$. Define range and source maps on $p\L$ by $r_p(\lambda) := \lambda(0,p)$, and $s_p(\lambda) := \lambda(d(\lambda) - p, d(\lambda))$ for all $\lambda \in p\L$, and define composition by $\lambda \circ_p \mu := \lambda \mu(p, d(\mu)) = \lambda(0, d(\lambda) - p) \mu$ whenever $s_p(\lambda) = r_p(\mu)$. Finally, define a degree map $d_p$ on $p\L$ by $d_p(\lambda) := d(\lambda) - p$ for all $\lambda \in p\L$.
\end{dfn}

\begin{prop} 
Let $(\Lambda,d)$ be a $k$-graph, and let $p\in\N^k$. Then $(p\Lambda,d_p)$ is a $k$-graph.
\end{prop}
\begin{proof} 
Define $\Obj(\Cc) := \L^p$, $\Mor(\Cc) := p\L$, $\cod_\Cc := r_p$, $\dom_\Cc := s_p$, $\Id_\Cc := \iota$, and $\circ_\Cc := \circ_p$. Then $\Cc = (\Obj(\Cc), \Mor(\Cc), \dom_\Cc, \cod_\Cc, \Id_\Cc, \circ_\Cc)$ is a category with morphisms $p\L$; it is straightforward to check that $\circ_p$ is associative using the factorisation property for $\L$. If $\lambda,\mu \in p\L$ and $s_p(\lambda) = r_p(\mu)$, then $d_p(\lambda\circ_p \mu) = d(\lambda \circ_p \mu) - p = d(\lambda\mu(p,d(\mu))) - p = (d(\lambda) + d(\mu) - p) - p = d(\lambda) - p + d(\mu) - p = d_p(\lambda) + d_p(\mu)$, and it follows that $d_p$ is a functor from $p\L$ to $\NN^k$. 

We need to check that the factorisation property holds for $p\Lambda$. Take any $\lambda\in p\Lambda$ and $m,n\in\N^k$ with $m+n=d_p(\lambda)$, so $d(\lambda)=m+p+n$. By the factorisation property for $\Lambda$ we have $\lambda = \lambda(0,m)\lambda(m,m+p)\lambda(m+p,m+p+n)$. But then $\lambda = \big(\lambda(0,m)\lambda(m,m+p)\big) \circ_p \big(\lambda(m,m+p)\lambda(m+p, m+p+n)\big)$ in $p\Lambda$, and $d_p(\lambda(0,m)\lambda(m,m+p))=m$ and $d_p(\lambda(m,m+p)\lambda(m+p, m+p+n))=n$. This decomposition is unique by the factorisation property for $\Lambda$.
\end{proof}

\begin{rmk}\label{rem:sourcesink}
If $\Lambda$ has no sources, then $p\L$ has no sources, and if $\L$ is row finite, then $p\Lambda$ is row finite.
\end{rmk}

\begin{prop}\label{prop:action} Let $(\Lambda,d)$ be a $k$-graph, and let $p,q\in\N^k$. Then $q(p\Lambda)=(q+p)\Lambda$.
\end{prop}
\begin{proof} 
By definition, we have $q(p\Lambda)^n=p\Lambda^{(n+q)}=\Lambda^{(n+q+p)}=(q+p)\Lambda^n$ for all $n \in \NN$.
Hence $q(p\Lambda)$ and $(q+p)\Lambda^n$ have identical elements. For the remainder of the proof, we write $s_q^{p\L}$, $r_q^{p\L}$, $\circ_q^{p\L}$, and $d_q^{p\L}$ for the source, range, composition and degree maps of the dual graph $q(p\L)$. 

Fix $\lambda \in \L^{n+p+q}$. We have that $s_{(q+p)}(\lambda) = \lambda(n, n+p+q)$ by definition, while $s^{p\L}_q(\lambda)$ is the final segment $\mu$ of $\lambda$ such that $d(\mu) - p = d_p(\mu) = q$; that is $d(\mu) = p+q$. Hence $s_{p+q}(\lambda) = s^{p\L}_q(\lambda)$. Similarly, $r_{p+q}(\lambda) = \lambda(0, p+q) = r^{p\L}_q(\lambda)$. Moreover, $d_{p+q}(\lambda) = d(\lambda) - (p+q) = d_p(\lambda) - q = d^{p\L}_q(\lambda)$. Since $\lambda$ was arbitrary, it follows that the range, source, and degree maps for $(p+q)\L$ and $q(p\L)$ agree.

This established, we have that $r_{p+q}(\lambda) = s_{p+q}(\mu)$ if and only if $r^{p\L}_q(\lambda) = s^{p\L}_q(\mu)$, in which case both $\lambda \circ_{p+q} \mu$ and $\lambda \circ^{p\L}_q \mu$ are equal to $\lambda\mu(p+q, d(\mu))$ by definition, completing the proof.
\end{proof}

\begin{theorem}\label{thm:cong}
Let $(\Lambda, d)$ be a row finite $k$-graph with no sources, and let $p \in \NN^k$. Let $\{s_\lambda : \lambda \in \L\}$ denote the universal generating Cuntz-Krieger $\L$-family in $C^*(\L)$, and let $\{t_\lambda : \lambda \in \L\}$ be the universal generating Cuntz-Krieger $p\L$-family in $C^*(p\L)$. For all $\lambda\in p\L$, define $r_\lambda := s_\lambda s_{s_{p}(\lambda)}^*$. There is an isomorphism $\phi : C^*(p\Lambda) \to C^*(\Lambda)$ such that $\phi(t_\lambda)=r_\lambda$ for all $\lambda \in p\L$.
\end{theorem}
\begin{proof}
First we show that the family $\{r_\lambda : \lambda\in p\L\}$ is a Cuntz-Krieger $p\L$-family. Since, for any $\beta\in p\L^0$, we have $s_\beta\ne 0$ it follows that $r_\beta = s_\beta s^*_\beta \ne 0$ and that it is a projection in $C^*(\Lambda)$. Furthermore, for distinct $\alpha,\beta\in p\L^0$, we have
\[
r_\alpha r_\beta = s_\alpha s_\alpha^* s_\beta s_\beta^*
=\delta_{\alpha\beta} s_\alpha s^*_\beta = \delta_{\alpha,\beta} r_{\alpha}.
\]
This establishes relation~(i).

For relation~(ii), let $\mu, \nu \in p\L$ with $r_{p}(\nu)=s_{p}(\mu)$, so $\mu\circ_{p}\nu = \mu\nu(p,d(\nu))$. Then,
\begin{equation} \label{eq:rearrange}
r_{\mu\circ_{p}\nu}
=s_{\mu\circ_{p}\nu}s^*_{s_{p}(\mu\circ_{p}\nu)}
=s_\mu s_{\nu(p,d(\nu))} s^*_{s_{p}(\nu)}
=s_\mu s^*_{s_{p}(\mu)} s_{s_{p}(\mu)} s_{\nu(p,d(\nu))} s^*_{s_{p}(\nu)}.
\end{equation}
But $s_{p}(\mu) = r_{p}(\nu) = \nu(0, p)$, so we can rewrite the right-hand side of~\eqref{eq:rearrange} to obtain $r_{\mu\circ_{p}\nu} = s_\mu s^*_{s_{p}(\mu)} s_{\nu} s^*_{s_{p}(\nu)} = r_\mu r_\nu$. This establishes relation~(ii).

Let $\lambda\in p\L$, say $d_{p}(\lambda) = n$. Then $r^*_\lambda r_\lambda = s_{s_{p}(\lambda)} s_\lambda^* s_\lambda s_{s_{p}(\lambda)}^* = s_{s_{p}(\lambda)} s^*_{s_{p}(\lambda)} = r_{s_{p}(\lambda)}$ by definition, establishing relation~(iii).

Finally, for relation~(iv), let $\beta \in p\L^0$ and let $n \in \NN^k$. Then
\[
r_\beta = s_\beta s_\beta^* 
= \sum_{\gamma\in s(\beta)\Lambda^n} s_\beta s_\gamma s_\gamma^* s_\beta^*
= \sum_{\lambda\in \beta\Lambda^n} s_\lambda s^*_\lambda.
\]
Applying the factorisation property and relation~(ii) for $C^*(\L)$ to the right-hand side then gives
\[
r_\beta = \sum_{\lambda\in \beta\Lambda^n} s_{\lambda(0,n)} s_{\lambda(n,n+p)} s^*_{\lambda(n,n+p)} s^*_{\lambda(0,n)},
\]
and then since each $s_{\lambda(n,n+p)} s^*_{\lambda(n,n+p)}$ is a projection, we obtain
\[
r_\beta
=\sum_{\lambda\in \beta\Lambda^n} (s_{\lambda(0,n)} s_{\lambda(n,n+p)} s^*_{\lambda(n,n+p)}) (s_{\lambda(n,n+p)} s^*_{\lambda(n,n+p)} s^*_{\lambda(0,n)}) 
=\sum_{\lambda\in \beta(p\Lambda^n)} r_\lambda r_\lambda^*,
\]
which establishes relation~(iv).

It follows from the universal property of $C^*(p\Lambda)$ that there exists a homomorphism $\phi : C^*(p\Lambda) \to C^*(\Lambda)$ satisfying $\phi(t_\lambda) = r_\lambda$ for all $\lambda \in p\L$. We claim that $\{r_\lambda : \lambda \in p\L\}$ generates $C^*(\Lambda)$. To see this, let $\sigma \in \L$ with $d(\sigma) = n$. An application of relation~(iv) for $C^*(\L)$ gives $s_\sigma = \sum_{\beta\in s(\sigma) \Lambda^{p}} s_\sigma s_\beta s_\beta^* 
= \sum_{\lambda \in \sigma\Lambda^{p}} s_\lambda s_{s_{p}(\lambda)}^*$, and this last is equal to $\sum_{\lambda \in \sigma\Lambda^{p}} r_\lambda$ by definition. Thus $\phi$ maps $C^*(p\Lambda)$ onto $C^*(\Lambda)$.

Now let $\gamma^\L$ denote the gauge action on $C^*(\L)$, and let $\gamma^{p\L}$ denote the gauge action on $C^*(p\L)$. For $z\in\mathbb{T}^k$ and $\lambda\in p\Lambda$, we have $\gamma^\L_z(r_\lambda) = \gamma^\L_z(s_\lambda s^*_{s_{p}(\lambda)}) = z^{d(\lambda)} s_\lambda (z^{d(s_{p}(\lambda))} s_{s_{p}(\lambda)})^* = z^{d(\lambda) - p} r_\lambda = z^{d_{p}(\lambda)} r_\lambda = \gamma^{p\L}(r_\lambda)$. Theorem~3.4 of \cite{KP} now establishes that $\phi$ is injective.
\end{proof}

\begin{rmk}
The hypotheses that $\L$ be row-finite and have no sources are crucial in Theorem~\ref{thm:cong}. To see why, notice that for $v \in \L^0$, the generator $s_v$ of $C^*(\L)$ is recovered in $C^*(p\L)$ as $\sum_{\beta \in p\L^0, r(\beta) = v} r_\beta$. However, even for 1-graphs, which contain sources or are not row-finite, the Cuntz-Krieger relations only insist that $p_v = \sum_{r(e) = v} s_e s^*_e$ when $r^{-1}(v)$ is finite and nonempty. 
\end{rmk}

\begin{lemma}\label{lem:overlap} 
Let $(\L,d)$ be a $k$-graph, and let $p \in \NN^k$. For each $n\in\N^k$ with $n\le p$ and $v,w\in p\Lambda^0$, there is at most one $\lambda\in v(p\Lambda^n)w$.
\end{lemma}
\begin{proof} 
Let $v,w \in p\L^0 = \L^p$ and suppose $\lambda \in v (p\L^n) w$. Then $\lambda \in \L^{n+p}$, $\lambda(0, p) = v$, and $\lambda(n, n+p) = w$. In particular, since $n \le p$ we have $\lambda(0,n) = \big(\lambda(0,p)\big)(0,n) = v(0,n)$, and then $\lambda = \lambda(0,n)\lambda(n,n+p) = v(0,n)w$, and hence is completely determined by $v$ and $w$.
\end{proof}

\begin{notation}
Let $(\L,d)$ be a $k$-graph. We write $M^\L_i$, $1 \le i \le k$ for the matrices in $M_{\L^0}(\NN)$ determined by $(M^\L_i)_{v,w} := |w \L^{e_i} v|$ for $w,v \in \L^0$, and we refer to these matrices as the \emph{coordinate matrices} of $\L$.
\end{notation}

\begin{cor} \label{cor:0-1}
Let $(\L,d)$ be a $k$-graph, and let and $p\in\N^k$ with $p_i\ge 1$ for $1 \le i \le k$. Then the coordinate matrices $M^{p\L}_i$ of $p\Lambda$ are $\{0,1\}$-matrices.
\end{cor}
\begin{proof}
The result is a direct consequence of Lemma~\ref{lem:overlap}.
\end{proof}

\section{$K$-theory}\label{sec:K-theory}

In this section we identify a class of $2$-graphs whose associated $C^*$-algebras are isomorphic to higher rank Cuntz-Krieger algebras in the sense of \cite{RS}, and use the results of \cite{RS} to calculate the $K$-theory of the $C^*$-algebras of such $2$-graphs. To state the main theorem for this section we employ the following notation: given square $n \times n$ matrices $M,N$, we write $\big[M\quad N\big]$ for the block $n \times 2n$ matrix whose first $n$ columns are those of $M$ and whose last $n$ columns are those of $N$. We also write $\One$ for the element $(1,1)$ of $\NN^2$. 

\begin{theorem}\label{thm:K-theory}
Let $(\L,d)$ be a $2$-graph which is finite and strongly connected as in Definition~\ref{dfn:structural conditions} and which has an aperiodic infinite path as in Definition~\ref{dfn:infinite paths}. Then $C^*(\L)$ is purely infinite, simple, unital and nuclear, and we have
\begin{align}
\rank(K_0(C^*(\L))) 
&= \rank(K_1(C^*(\L))) \notag \\ 
&= \rank\big(\coker\big[I - M^{\One\L}_1\quad I - M^{\One\L}_2\big]\big) \label{eq:K-th1} \\
& \qquad\qquad + \rank\big(\coker\big[I - \big(M^{\One\L}_1\big)^t\quad I - \big(M^{\One\L}_2\big)^t\big]\big); \notag \displaybreak[0]\\
\tor(K_0(C^*(\L))) &\cong \tor\big(\coker\big[I - M^{\One\L}_1\quad I - M^{\One\L}_2\big]\big);\text{ and} \\
\tor(K_1(C^*(\L))) &\cong \tor\big(\coker\big[I - (M^{\One\L}_1\big)^t\quad I - \big(M^{\One\L}_2\big)^t\big]\big) \label{eq:K-th4}.
\end{align}
\end{theorem}

The remainder of this section constitutes the proof of Theorem~\ref{thm:K-theory}. We begin by recalling some definitions from \cite{RS}. Let $A$ be a finite set, and let $M_1, M_2$ be $A \times A$ matrices with entries in $\{0,1\}$. For $n \in \NN^k$, let $W_n := \{w : [0,n] \to A : M_j(w(l+e_j), w(l)) = 1$ whenever $l, l+e_j \in [0,n]\}$; we refer to the elements of $W_n$ as \emph{allowable words of shape $n$}, and write $W$ for the collection $\bigcup_{n \in \NN^k} W_n$ of all \emph{allowable words}. For $u \in W$, write $S(u)$ for the shape of $u$; that is, $S(u)$ is the unique element of $\NN^k$ such that $u \in W_{S(u)}$.  Notice that $W_0$ is just $A$. The matrices $M_1, M_2$ are said to satisfy~(H0)--(H3) if
\begin{itemize}
\item[(H0)] Each $M_i$ is nonzero;
\item[(H1a)] $M_1 M_2 = M_2 M_1$;
\item[(H1b)] $M_1 M_2$ is a $\{0,1\}$-matrix;
\item[(H2)] the directed graph with a vertex for each $a \in A$ and a directed edge $(a,i,b)$ from $a$ to $b$ for each $a,i,b$ such that $M_i(b,a) = 1$, is irreducible; and
\item[(H3)] for each $m \in \ZZ^2 \setminus \{0\}$, there exists a word $w \in W$ and elements $l_1,l_2$ of $\NN^2$ with $0 \le l_1,l_2 \le S(w)$ such that $l_2 - l_1 = m$ and $w(l_1) \not= w(l_2)$. 
\end{itemize}

\begin{notation}
If $(\L,d)$ is a $2$-graph such that the coordinate martices $M^\L_1$ and $M^\L_2$ are $\{0,1\}$-matrices, we write $W^\L_n$ and $W^\L$ for the collection of allowable words of shape $n$ and for the collection of all allowable words respectively. For $\lambda \in \L$, let $w^\L_\lambda$ be the word in $W^{\L}_{d(\lambda)}$ given by $w^\L_\lambda(m) = s(\lambda(0,m))$ for $0 \le m \le d(\lambda)$. Since each $M^{\L}_i$ is a $\{0,1\}$-matrix, the map $\lambda \mapsto w^\L_\lambda$ is a bijection between $\L^n$ and $W^{\L}_n$ for all $n \in \NN^k$.
\end{notation}

\begin{prop} \label{prp:RS algebras}
Let $(\L,d)$ be a finite $2$-graph with no sources, and let $M_1^{\One\L}$ and $M_2^{\One\L}$ be the matrices associated to the higher-edge graph $\One\L$. Then
\begin{itemize}
\item[(1)] $M^{\One\L}_1, M^{\One\L}_2$ satisfy \rm{(H0)}, \rm{(H1a)}, and \rm{(H1b)};
\item[(2)] $M^{\One\L}_1, M^{\One\L}_2$ satisfy \rm{(H2)} if and only if $\L$ is strongly connected; and
\item[(3)] if $M^{\One\L}_1, M^{\One\L}_2$ satisfy \rm{(H2)}, then they satisfy \rm{(H3)} if and only if $\L$ has an aperiodic infinite path.
\end{itemize}
\end{prop}
\begin{proof}
For (1), note that each $M^{\One\L}_i$ is a finite square matrix over $\One\L^0$ by definition, and has entries in $\{0,1\}$ by Corollary~\ref{cor:0-1}. It is easy to see that $(M_i^{\One\L} M_{3-i}^{\One\L})_{v,w} = |\{(\alpha,\beta) \in w(\One\L^{e_{3-i}}) \times (\One\L^{e_i})v : r(\alpha) = s(\beta)\}| = |w (\One\L^\One) v|$ for $i= 1,2$ and this establishes~(H1a). The same calculation combined with Lemma~\ref{lem:overlap} establishes~(H1b).

For (2), notice that $M^{\One\L}_1, M^{\One\L}_2$ satisfy~(H2) if and only if for every $v,w \in \One\L^0$ there exist elements $\alpha_1, \dots, \alpha_k$ in $\One\L^{(1,0)} \cup \One\L^{(0,1)}$ such that $r(\alpha_1) = v$, $s(\alpha_k) = w$, and $r(\alpha_{i+1}) = s(\alpha_i)$ for $1 \le i \le k-1$. 

So suppose first that $M^{\One\L}_1, M^{\One\L}_2$ satisfy~(H2), and let $v,w \in \L^0$. Since $\L$ has no sources, there exist $\mu,\nu \in \L^\One$ with $r(\mu) = v$ and $r(\nu) = w$; so $\mu,\nu \in \One\L^0$ by definition, and~(H2) ensures that there is a path $\alpha_1, \dots \alpha_k$ from $\mu$ to $\nu$ in $\One\L^{(1,0)} \cup \One\L^{(0,1)}$. By definition of $\One\L$, the path $\alpha_1\dots\alpha_k$ in $\One\L$ is a path $\lambda \in \L$ with $d(\lambda) = d_\One(\alpha_1\dots\alpha_k) + \One$, and such that $\lambda(0,\One) = \mu$ and $\lambda(d(\lambda) - \One, d(\lambda)) = \nu$. But then $\lambda(0, d(\lambda)-\One) \in v\L w$. Since $v,w \in \L^0$ were arbitrary, it follows that $\L$ is strongly connected.

Now suppose that $\L$ is strongly connected, and fix $\mu,\nu \in \One\L^0$. Since $\L$ is strongly connected, there is a path $\lambda \in s(\mu) \L r(\nu)$, and then $\tau := \mu\lambda\nu$ belongs to $\mu (\One\L) \nu$ with $d_\One(\mu\lambda\nu) = d(\lambda) + \One$. Any factorisation of $\tau$ into segments from $\One\L^{(1,0)} \cup \One\L^{(0,1)}$ now gives a path in $\One\L^{(1,0)} \cup \One\L^{(0,1)}$ from $\nu$ to $\mu$, so $M^{\One\L}_1, M^{\One\L}_2$ satisfy~(H2).

Finally, for (3), assume that $M^{\One\L}_1, M^{\One\L}_2$ satisfy \rm{(H2)}, so $\L$ is strongly connected by part~(2). For $x \in \L^\infty$, define $\One x \in \One\L^\infty$ by $(\One x)(m,n) := x(m, n+\One)$. It is easy to see that the map $x \mapsto \One x$ is a bijection between $\L^\infty$ and $\One\L^\infty$. 

Claim: $x \in \L^\infty$ is aperiodic if and only if $\One x \in \One\L^\infty$ is aperiodic. To see this, let $m,n \in \NN^k$, and fix $x \in \L^\infty$. By definition, we have 
\begin{align}
\sigma^m(\One x) = \sigma^n(\One x)
&\iff (\One x)(s + m, t + m) = (\One x)(s + n, t+ n)\quad\text{for $s \le t$} \notag\\
&\iff x(s + m, t + m + \One) = x(s + n, t + n + \One)\quad\text{for $s \le t$} \label{eq:px periodic} 
\end{align}
Now if $x(s+m, t+m+\One) = x(s+n, t+n+\One)$ for all $s \le t \in \NN^2$, then the uniqueness of factorisations in $\L$ ensures that $x(s+m, t+m) = x(s+n, t+n)$ for all $s \le t \in \NN^2$. Conversely if $x(s+m, t+m) = x(s+n, t+n)$ for all $s \le t \in \NN^2$, then replacing $t$ with $t+\One$ gives $x(s+m, t+m+\One) = x(s+n, t+n+\One)$ for all $s \le t \in \NN^2$. Hence~\eqref{eq:px periodic} shows that
\begin{align*}
\sigma^m(\One x) = \sigma^n(\One x)
&\iff x(s+m, t+m) = x(s+n, t+n)\quad\text{for $s \le t \in \NN^2$} \\
&\iff \sigma^m(x) = \sigma^n(x),
\end{align*}
establishing the claim. Thus it suffices to show that $M_i^{\One\L}$ satisfy~(H3) if and only if $\One\L^\infty$ has an aperiodic element.

Suppose first that there exists an aperiodic path $x \in \One\L^\infty$. Fix $m \in \ZZ^2$, and write $m = m_+ - m_-$ where $m_+, m_- \in \NN^2$. Since $|v(\One\L^{e_i})w| \in \{0,1\}$ for all $v,w \in \One\L^0$, $i = 1,2$, we have that $x$ is completely determined by its restriction to the objects of $\Omega_2$; that is, by the function from $\NN^2$ to $\L^0$ given by $n \mapsto x(n)$. Since $x$ is aperiodic, it follows that $\sigma^{m_+}(x)(n) \not= \sigma^{m_-}(x)(n)$ for some $n \in \NN^2$. But then with $N := n + m_-$, we have $x(N + m_+ - m_-) \not = x(N)$, and $w := x|_{[0, N + m_+ - m_-]} \in W^{\One\L}_{N+ m_+ - m_-}$ satisfies $w(N) \not= w(N + m)$. Since $m \in \ZZ^2$ was arbitrary, this establishes that $M^{\One\L}_1, M^{\One\L}_2$ satisfy \rm{(H3)}.

Now suppose that $M^{\One\L}_1, M^{\One\L}_2$ satisfy \rm{(H3)}. For each $m \in \ZZ^2\setminus\{0\}$, fix $w_m \in W^{\One\L}$ and $l_m \in \NN^2$ such that $0 \le l_m, l_m + m \le S(w_m)$ and $w_m(l_m) \not= w_m(l_m + m)$. Let $\lambda_m$ be the unique path in $\One\L$ such that $w_m = w^{\One\L}_{\lambda_m}$. We will construct an infinite path $x$ which contains infinitely many occurrences of each $\lambda_m$; this will ensure that there is no $m$ for which a sufficiently large shift of $x$ has period $m$, and hence that $x$ is aperiodic. The details of this construction, and the verification that the resulting $x$ is aperiodic constitute the remainder of the proof.

Let $\{m_i : i \in \NN\}$ be a listing of $\ZZ^2 \setminus \{0\}$. Fix an arbitrary $v \in \One\L^0$, and for each $i \in \NN$, let $\alpha_i$ be any element of $v (\One\L) r(\lambda_{m_i})$, and let $\beta_i$ be any element of $s(\lambda_{m_i})(\One\L) v$ with the property that $d_\One(\alpha_i \lambda_{m_i} \beta_i) \ge \One$; this is possible because $\L$ is strongly connected and has no sources.

For $i \in \NN$, let $\rho_i := \alpha_i\lambda_{m_i}\beta_i$, and let $\tau_i := \rho_1 \rho_2 \dots \rho_i$. Let $x$ be the infinite path $x := \tau_1\tau_2\tau_3\cdots$. We claim that $x$ is aperiodic.

To see this, let $s,t \in \NN^2$, and let $I_{s,t}$ be the element of $\NN$ such that $m_{I_{s,t}} = t-s$. Let $J := \max\{s_1,s_2,t_1,t_2\}$; since $d_\One(\rho_i) \ge (1,1)$, we have that $i \ge J$ implies $d_\One(\tau_1\cdots\tau_i) \ge s,t$. Let $K := \max\{I_{s,t}, J + 1\}$, and define $N := d_\One(\tau_1 \cdots \tau_{K-1}) + d_\One(\rho_1 \cdots \rho_{I_{s,t}-1}) + d(\alpha_{I_{s,t}}) + l_{t-s} - s$. We have $N \ge 0$ by choice of $K$, and
\begin{align*}
\sigma^s(x)(N)
&= x(N + s) \\
&= x(d_\One(\tau_1 \cdots \tau_{K-1}) + d_\One(\rho_1 \cdots \rho_{I_{s,t}-1}) + d(\alpha_{I_{s,t}}) + l_{t-s}) \\
&= \lambda_{m_{I_{s,t}}}(l_{t-s}).
\end{align*}
A similar calculation shows that $\sigma^t(x)(N) = \lambda_{m_{I_{s,t}}}(l_{t-s} + (t-s))$, and hence $\sigma^s(x)(N) \not = \sigma^t(x)(N)$ by our chice of $\lambda_{m_{I_{s,t}}}$. It follows that $\sigma^s(x) \not= \sigma^t(x)$, and since $s,t \in \NN^2$ were arbitrary, that $x$ is aperiodic.
\end{proof}

\begin{notation}
Let $\L$ be a finite strongly connected $2$-graph with an aperiodic infinite path. We write $\mathcal{A}^{\One\L}$ for the $C^*$-algebra associated to $M^{\One\L}_i$ as in \cite{RS}. That is, $\mathcal{A}^{\One\L}$ is the universal $C^*$-algebra generated by a family $\{s_{u,v} : u,v \in W^{\One\L}, u(S(u)) = v(S(v))\}$ of partial isometries satisfying
\begin{align}
s_{u,v} &= s^*_{v,u}\quad\text{for $u,v \in W^{\One\L}$;} \label{eq:A1}\\
s_{u,v} s_{v,w} &= s_{u,w}\quad\text{for $u,v,w \in W^{\One\L}$;} \\
s_{u,v} &= \sum_{w \in W^{\One\L}_{e_j}, u(S(u)) = w(0)} s_{uw} s^*_{vw}\quad\text{for $u,v, \in W^{\One\L}$, $j \in \{1,2\}$; and}  \\
s_{a,a} s_{b,b} &=0 \quad\text{for distinct $a,b \in W^{\One\L}_0$.} \label{eq:A4}\displaybreak[0]
\end{align}
\end{notation}

\begin{lemma}
Let $(\L,d)$ be a finite strongly connected $2$-graph which has an aperiodic infinite path. Then $C^*(\L)$ is isomorphic to $\mathcal{A}^{\One\L}$.
\end{lemma}
\begin{proof}
The factorisation property ensures that if $\L$ is strongly connected and contains an infinite path, then $\L$ has no sources. By Theorem~\ref{thm:cong}, we have that $C^*(\L)$ is isomorphic to $C^*(\One\L)$, so it suffices to show that $C^*(\One\L)$ is isomorphic to $\mathcal{A}^{\One\L}$. It is easy to check using Definition~\ref{dfn:CKfam}(i)--(iv), relations \eqref{eq:A1}--\eqref{eq:A4}, and the universal properties of $\mathcal{A}^{\One\L}$ and $C^*(\One\L)$ that there exists a homomorphism $\pi : \mathcal{A}^{\One\L} \to C^*(\One\L)$ satisfying $\pi(s_{w^{\One\L}_\lambda,w^{\One\L}_\mu}) = s_\lambda s^*_\mu$ for all $\lambda,\mu \in \One\L$, and that there exists a homomorphism $\psi : C^*(\One\L) \to \mathcal{A}^{\One\L}$ satisfying $\psi(s_\lambda) := s_{w^{\One\L}_\lambda, w^{\One\L}_{s(\lambda)}}$. Since these two homomorphisms are mutually inverse, the result follows.
\end{proof}

\begin{rmk}
The argument of statement~(2) of Proposition~\ref{prp:RS algebras} shows that if $\L$ has no sources, then for any $q \ge \One$, the coordinate matrices of $q\L$ will satisfy~(H2) only if $\L$ is strongly connected and has no sources. In particular, there exists $q \in \NN^2$ such that $M^{q\L}_i$ satisfy~(H0)--(H3) if and only if $M^{\One\L}_i$ satisfy~(H0)--(H3).
\end{rmk}

\begin{proof}[Proof of Theorem~\ref{thm:K-theory}]
Theorem~5.9, Proposition~5.11, and Corollary~6.4 of \cite{RS'} combined with the previous two results show that $C^*(\L)$ is simple, purely infinite and nuclear. We have that $C^*(\L)$ is unital with $1_{C^*(\L)} = \sum_{v \in \L^0} s_v$. Proposition~2.14 of \cite{RS} establishes \eqref{eq:K-th1}--\eqref{eq:K-th4}. 
\end{proof}

\begin{rmks}

(1) The proof of \cite[Proposition~2.14]{RS} does not make any use of relations (H2)~and~(H3). Hence the formulas for $K_*(C^*(\L))$ in Theorem~\ref{thm:K-theory} hold when $\L$ is a finite $k$-graph with no sinks or sources, even if it is not strongly connected and does not have an aperiodic infinite path. However, in this case $C^*(\L)$ is not necessarily simple and purely infinite, and so is not determined up to isomorphism by its $K$-theory.

(2) The formulas for $K_*(C^*(\L))$ given in Theorem~\ref{thm:K-theory} are in terms of the coordinate matrices $M^{\One\L}_i$ of the dual $k$-graph. Proposition~5.1 of \cite{Gw} shows that the same formulas hold if all instances $M_i^{\One\L}$ are replaced with $M^\L_i$, but it is unclear how to show this directly.
\end{rmks}


\begin{thebibliography}{00}

\bibitem{B} T. Bates, \emph{Applications of the gauge-invariant uniqueness theorem for graph algebras}, Bull. Austral. Math. Soc. {\bf 65} (2002), 55--67.

\bibitem{BP} T. Bates and D. Pask, \emph{Flow equivalence of graph algebras}, Ergodic Theory Dynam. Systems {\bf 24} (2004), in press.

\bibitem{BPRS} T. Bates, D. Pask, I. Raeburn, and W. Szyma\'nski, \emph{The $C^*$-algebras of row--finite graphs}, New
York J. Math.  {\bf 6} (2000), 307--324.

\bibitem{CK} J. Cuntz and W. Krieger, \emph{A class of $C^*$-algebras and topological Markov chains,} Invent.  Math. 
{\bf 56} (1980), 251--268.

\bibitem{EW} M. Enomoto and Y. Watatani, \emph{A graph theory for $C^* $-algebras}, Math.  Japon.  {\bf 25} (1980),
435--442.

\bibitem{Gw} D. G. Evans, \emph{On higher-rank graph $C^*$-algebras}, Ph.D. Thesis, Univ. Wales, 2002.

\bibitem{FLR} N. J. Fowler, M. Laca, and I. Raeburn, {\em The $C^*$-algebras of infinite graphs}, Proc.  Amer.  Math. 
Soc.  {\bf 128} (2000), 2319--2327.

\bibitem{KP} A. Kumjian and D. Pask, \emph{Higher rank graph $C^*$-algebras}, New York J. Math., {\bf 6} (2000) 1--20.

\bibitem{KPR} A. Kumjian, D. Pask, and I. Raeburn, {\em Cuntz-Krieger algebras of directed graphs}, Pacific J. Math. {\bf 184} (1998), 161--174.

\bibitem{KPRR} A. Kumjian, D. Pask, I. Raeburn, and J. Renault, \emph{Graphs, groupoids and Cuntz-Krieger algebras}, J. Funct.  Anal.  {\bf 144} (1997), 505--541.

\bibitem{MRS} M.H. Mann, I. Raeburn, and C.E. Sutherland, \emph{Representations of finite groups and Cuntz-Krieger algebras}, Bull. Austral. Math. Soc. {\bf 46} (1992), 225--243.

\bibitem{P} N.C. Phillips, \emph{A classification theorem for nuclear purely infinite simple $C^*$-algebras,}
Documenta Math. {\bf 5} (2000), 49--114.

\bibitem{RaSz} I. Raeburn and W. Szyma\'nski, \emph{Cuntz-Krieger algebras of infinite graphs and matrices}, Trans. 
Amer. Math. Soc. {\bf 356} (2004), 39--59.

\bibitem{RSY1} I. Raeburn, A. Sims, and T. Yeend, \emph{Higher-rank graphs and their $C^*$-algebras}, Proc. Edinb. Math. Soc {\bf 46} (2003), 99--115.

\bibitem{RS'} G. Robertson and T. Steger, \emph{Affine buildings, tiling systems and higher rank Cuntz-Krieger algebras}, J. reine angew. Math. {\bf 513} (1999), 115--144.

\bibitem{RS} G. Robertson and T. Steger, \emph{Asyptotic $K$-theory for groups acting on $\tilde A_2$ bulidings}, Can. J. Math. {\bf 53} (2001), 809--833.

\bibitem{Sz} W. Szyma\'nski, \emph{The range of $K$-invarants for $C^*$-algebras of infinite graphs}, Indiana Univ. 
Math.  J. {\bf 51} (2002), 239--249.

\end{thebibliography}
\end{document}